\documentclass{elsarticle}
\input{amssym.def}
\input{amssym.tex}
\title{Reinforcement of a plate weakened by multiple holes with several patches for different types of a plate-patch attachment}

\author{A. Y. ZEMLYANOVA}
\address{Department of Mathematics, IAMCS, Texas A\&M University\\
Mailstop 3368, Texas A\&M University, College Station, TX 77843-3368 USA\\
azem@math.tamu.edu}

\date{}

\newcommand{\I}{\mathop{\rm Im}\nolimits}
\newcommand{\R}{\mathop{\rm Re}\nolimits}

\newcommand{\beqa}{\begin{eqnarray}}
\newcommand{\eeqa}[1]{\label{#1}\end{eqnarray}}
\newcommand{\bequ}{\begin{equation}}
\newcommand{\eequ}[1]{\label{#1}\end{equation}}

\newcommand{\beq}{\begin{equation}}
\newcommand{\eeq}{\end{equation}}
\newcommand{\barr}{\begin{eqnarray}}
\newcommand{\earr}{\end{eqnarray}}
\newcommand{\beqn}{\begin{equation*}}
\newcommand{\eeqn}{\end{equation*}}
\newcommand{\barrn}{\begin{eqnarray*}}
\newcommand{\earrn}{\end{eqnarray*}}

\begin{document}

\noindent

\begin{abstract}

The most general situation of the reinforcement of a plate with multiple holes by several patches is considered. There is no restriction on the number and the location of the patches. Two types of the patch attachment are considered: only along the boundary of the patch or both along the boundary of the patch and the boundaries of the holes which this patch covers. The unattached boundaries of the holes may be loaded with given in-plane stresses. The mechanical problem is reduced to the system of singular integral equations which can be further reduced to the system of Fredholm equations. A new numerical procedure for the solution of the system of singular integral equations is proposed in this paper. It can be observed that this procedure provides a noticeable improvement from the one presented in Zemlyanova \cite{Zem2007}, Zemlyanova and Silvestrov \cite{ZemSil2007}, \cite{ZemSil2009} and allows to achieve a significantly better numerical convergence with less computational effort. 

\end{abstract}


\maketitle

\pagestyle{myheadings}
\thispagestyle{plain}
\markboth{A. Y. ZEMLYANOVA}{Reinforcement of a plate weakened by multiple holes with several patches}

\setcounter{equation}{0}

\section{Introduction}

Patch repair is one of the most common techniques used to strengthen thin constructions containing imperfections such as cracks and holes. There is a significant amount of scientific literature dedicated to the topic with different types of the plate-patch attachment being considered. In particular,  the plate and the patches can be attached to each other using adhesive applied along parts of their surfaces (two-dimensional) or along some contours (one-dimensional), or only at certain point by using rivets. The detailed review of the literature on the bonded (two-dimensional) repair of the cracked plates can be seen in \cite{BakJon1988}, \cite{ErdArin1972}, \cite{Rose1988} and many others. Reinforcement of plates with holes is studied less extensively. Engels et al. \cite{EngZakBec2001} have considered a bonded repair of an anisotropic plate with a hole. Some problems on the reinforcement of a cracked plate with circular holes have been studied numerically in \cite{MWCh1975}. The patch repair of a single hole in a thin plate reinforced by a single patch is considered for two different types of the plate-patch attachment in  \cite{Zem2007}, \cite{ZemSil2007}. To the author's knowledge, a patch repair of multiple holes in a plate by multiple patches of arbitrary shapes has not been a subject of many studies with an exception of \cite{ZemSil2009}.

In this paper, a reinforcement of a plate with multiple holes with several patches is considered for two types of the plate-patch attachment: only along the boundary of the patch or along both the boundary of the patch and the boundaries of the holes covered by the patch. It is assumed that the patch is loaded with in-plane stresses applied to the free unattached boundaries of the holes and at infinity point of the plate. The location, the quantity and the shape of the patches can be arbitrary with the only condition that the boundaries of the patches and the holes are smooth and do not touch each other or intersect.

The complex analysis methods are used to solve the problem. First, the Muskhelishvili's formulas \cite{Mus1963} are used to describe the stresses and the displacements in the plate and the patches through the complex potentials in the plate and the patches. The singular integral representations of the complex potentials are used then to reformulate the problem as a system of singular integral equations. These representations were first proposed by Savruk \cite{Savruk1981} to study the problems for cracked plates. The uniqueness of the solution of the resulting system of the singular integral equations can be proved in a manner similar to \cite{Zem2007}, \cite{ZemSil2007}, \cite{ZemSil2009}.

One of the main goals of this paper is to propose an alternative technique for a numerical solution of the derived system of singular integral equations. The method of mechanical quadratures applied for the solution of similar systems in \cite{Zem2007}, \cite{ZemSil2007}, \cite{ZemSil2009} is difficult to implement for a large number of contours due to decreasing precision. The idea behind the method proposed here is in the approximation of the unknown functions by trigonometric polynomials. Substituting these approximations into the system of singular integral equations reduces this system to the system of linear algebraic equations for the coefficients of the trigonometric polynomials. Numerical examples show that even the systems of the relatively small order are sufficient to achieve a good accuracy of the results. By considering numerical examples it can be shown that this method demonstrates a noticeably better convergence, larger independence from the symmetry of the construction and, hence, can be applied to a wider range of practical problems.

\setcounter{equation}{0}

\section{Boundary conditions}
 
\begin{figure}
	\centering
		\scalebox{0.5}{\includegraphics{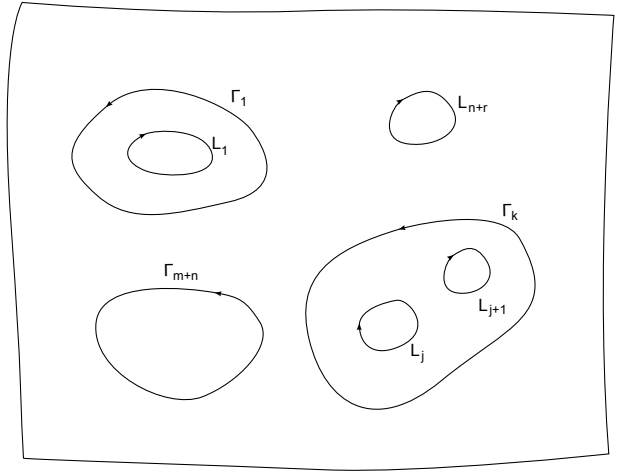}}
	\caption{A plate $S$ with holes and attached patches $S_k$.}
	\label{fig1}
\end{figure}

Consider an infinite thin elastic plate $S$ weakened by several holes with the boundaries $L_j$, $j=1,\ldots,n+r$. Thin elastic patches $S_k$, $k=1,\ldots, n+m$, with the boundaries $\Gamma_k$ are used to reinforce the plate $S$ (fig. \ref{fig1}). It is assumed that the patch $S_k$, $k=1,\ldots,n$, covers the hole $L_k$ and is attached to the plate $S$ both along the boundary of the hole $L_k$ and the boundary of the patch $\Gamma_k$. Other patches $S_k$, $k=n+1,\ldots,n+m$, are attached to the plate only along their boundaries $\Gamma_k$ and can be located anywhere in the plate $S$, in particular, they may cover some of the holes $L_j$, $j=n+1,\ldots,n+r$. Given in-plane stresses $p_j(t)$ act on the free boundaries of the holes $L_j$, $j=n+1,\ldots,n+r$. It is assumed that the lines $L_j$ and $\Gamma_k$ are smooth Lyapunov curves and do not touch or intersect each other.
 
The plate and the patches are homogeneous and isotropic. Their thicknesses, shear moduli and Poisson ratios are given by $h$, $\mu$, $\nu$ and $h_k$, $\mu_k$, $\nu_k$, $k=1,\ldots,m+n$, correspondingly. The principal in-plane stresses $\sigma_1^{\infty}$ and $\sigma_2^{\infty}$ are applied at infinity of the plate and act in the directions constituting the angles $\alpha$ and $\alpha+\pi/2$ with the positive direction of the real axis.

Assume that the patches and the plate are joined perfectly along the junction lines. This means that the displacements are equal both in the patch and the plate from both sides of the junction lines and the stress equilibrium conditions are satisfied on the junction lines. These conditions together with the given stresses on the free boundaries of the holes lead to the following set of boundary conditions:

\begin{equation}
\label{2.1}
(u+iv)^{+}(t)=(u+iv)_k^{+}(t)=(u+iv)_k^{-}(t), \, \, t\in L_k, \,\, k=1,\ldots,n,
\end{equation}    
\begin{equation}
\label{2.2}
h(\sigma_n+i\tau_n)^{+}(t)+h_k(\sigma_n+i\tau_n)^{+}_k(t)=h_k(\sigma_n+i\tau_n)^{-}_k(t),
\end{equation}  
$$
t\in L_k, \,\, k=1,\ldots,n,
$$
\begin{equation}
\label{2.3}
(\sigma_n+i\tau_n)^{+}_k(t)=p_k(t), \,\, t\in L_k, \,\, k=n+1,\ldots,n+r,
\end{equation}  
\begin{equation}
\label{2.4}
(u+iv)^{+}(t)=(u+iv)^{-}(t)=(u+iv)_k^{+}(t), \, \, t\in \Gamma_k, \,\, k=1,\ldots,n+m,
\end{equation}    
\begin{equation}
\label{2.5}
h(\sigma_n+i\tau_n)^{+}(t)+h_k(\sigma_n+i\tau_n)^{+}_k(t)=h(\sigma_n+i\tau_n)^{-}(t),
\end{equation}
$$
t\in \Gamma_k, \,\, k=1,\ldots,n+m,
$$
where $(u+iv)(t)$ is the vector of the displacements at a point $t$ of the plate $S$ or the patch $S_k$, and $\sigma_n$ and $\tau_n$ are the tensile and the shear components of the stress vector acting on the tangent line to the curves $L_k$ or $\Gamma_k$. Parameters without a subscript are related to the plate $S$; parameters with a subscript ``$k$" are related to the patch $S_k$. Here and henceforth, superscripts ``$+$" and ``$-$" denote the limit values of the stresses, the displacements and other parameters from the left-hand and the right-hand sides correspondingly of the lines $L_k$ or $\Gamma_k$. The direction of the curve $L_k$ is chosen to be clockwise and of the curve $\Gamma_k$ counterclockwise (fig. 1).

\section{Integral representations of complex potentials}
The stresses and the derivatives of the displacements in the plate $S$ or in the patch $S_k$ ($k=1,\ldots,m+n$) on any contour $L'$ lying in $S$ or in $S_k$ correspondingly can be obtained from the Muskhelishvili's complex potentials $\Phi(z)$, $\Psi(z)$ by the following formulas \cite{Mus1963}:\\
\begin{equation}
\label{3.1}
(\sigma_n+i\tau_n)(t)=\Phi(t)+\overline{\Phi(t)}+\frac{\overline{dt}}{dt}(t\overline{\Phi'(t)}+\overline{\Psi(t)}),
\end{equation}
$$
2\mu\frac{d}{dt}(u+iv)(t)=\kappa\Phi(t)-\overline{\Phi(t)}-\frac{\overline{dt}}{dt}(t\overline{\Phi'(t)}+\overline{\Psi(t)}),\,\,\,\,t\in L'.
$$
For the patch $S_k$, all the parameters and the functions in these formulas should be taken with the subscript ``$k$". In the case of plane stress take $\kappa=(3-\nu)/(1+\nu)$ and $\kappa_k=(3-\nu_k)/(1+\nu_k)$, $k=1,\ldots,m+n$.

We will derive the integral representations of the complex potentials in the spirit of \cite{Savruk1981}. The main idea is to divide the ``plate-patches" system into separate components and introduce the unknown functions describing the jumps of the stresses or the derivatives of the displacements on each of the lines $L_j$ and $\Gamma_k$. These unknown functions may be chosen in such a manner that some of the conditions (\ref{2.1}) - (\ref{2.5}) are satisfied automatically. 

First, consider the infinite plate $S$ with the holes $L_j$ which is subjected to the given in-plane stresses at infinity and on the lines $L_j$, $j=n+1,\ldots,n+r$. The stresses in the plate $S$ have a jump discontinuity on the curve $\gamma=\cup_{k=1}^{n+m} \Gamma_k$ lying in the interior of the plate, and the displacements are continuous on these curves:
\begin{equation}
\label{3.2}
q(t)=\left((\sigma_n+i\tau_n)^+(t)-(\sigma_n+i\tau_n)^-(t)\right)/2, \,\, t\in \gamma.
\end{equation}
Extend formally the plate $S$ to the full complex plane so that the stresses are continuous through the curve $L=\cup_{j=1}^{n+r}L_j$ and the derivatives of the displacements have a jump discontinuity along this line:
\begin{equation}
\label{3.3}
g'(t)=\frac{2\mu}{i(\kappa+1)}\frac{d}{dt}\left((u+iv)^+(t)-(u+iv)^-(t)\right), \,\, t\in L.
\end{equation}
Then the complex potentials $\Phi(z)$ and $\Psi(z)$ can be taken in the form \cite{Savruk1981}:\\
$$
\Phi(z)=\Gamma-\sum_{k=n+1}^{n+r}\frac{Q_k}{z-z_k}+\frac{1}{2\pi}\int_L\frac{g'(t)dt}{t-z}+\frac{(\kappa+1)^{-1}}{\pi i}\int_{\gamma}\frac{q(t)dt}{t-z},
$$
\begin{equation}
\label{3.4}
\Psi(z)=\Gamma'+\sum_{k=n+1}^{n+r}\frac{\kappa\bar{Q}_k}{z-z_k}+\frac{1}{2\pi}\int_L\left( \frac{\overline{g'(t)dt}}{t-z}-\frac{\bar{t}g'(t)dt}{(t-z)^2}\right)
\end{equation}
$$
{}+\frac{(\kappa+1)^{-1}}{\pi i}\int_{\gamma}\left(\frac{\kappa\overline{q(t)dt}}{t-z}-\frac{\bar{t}q(t)dt}{(t-z)^2}\right), \,\, z\in S,
$$
$$
\Gamma=(\sigma_1^{\infty}+\sigma_2^{\infty})/4, \,\,\,\, \Gamma'=(\sigma_2^{\infty}-\sigma_1^{\infty})e^{-2i\alpha}/2,
$$
$$
Q_k=\frac{X_k+iY_k}{2\pi(1+\kappa)}, \,\,\,\, X_k+iY_k=-i\int_{L_k}p_k(t)dt,
$$
where $z_k$ ($k=n+1,\ldots, n+r$) is the arbitrarily fixed point inside of the contour $L_k$.

Formally extend each patch $S_k$, $k=1,\ldots,n$, to the full complex plane so that outside of the patch all the displacements and the stresses are equal to zero. The displacements in the patch $S_k$ are continuous through the line $L_k \subset S_k$, and the stresses have a jump discontinuity on this line:
\begin{equation}
\label{3.5}
q_k(t)=\left((\sigma_n+i\tau_n)^+_k(t)-(\sigma_n+i\tau_n)^-_k(t)\right)/2, \,\, t\in L_k.
\end{equation}
In the complex plane thus obtained from the patch $S_k$ the jump of the vector of the stresses $(\sigma_n+i\tau_n)_k(t)$ acting on the tangent line to the contour $\Gamma_k$ is equal to $(\sigma_n+i\tau_n)^+_k(t)$, and the jump of the vector of the displacements $(u+iv)(t)$ is equal to $(u+iv)^+_k(t)$:
\begin{equation}
\label{3.6}
g'_k(t)=\frac{2\mu_k}{i(\kappa_k+1)}\frac{d}{dt}(u+iv)^+_k(t), \,\, t\in \Gamma_k.
\end{equation}
According to the condition (\ref{2.5}) and the formula (\ref{3.2}) we have $(\sigma_n+i\tau_n)^+_k(t)=-2d_k^{-1}q(t)$, $t\in \Gamma_k$, $d_k=h_k/h$. Therefore, the complex potentials $\Phi_k(z)$, $\Psi_k(z)$ can be taken in the form \cite{Savruk1981}:\\
\begin{equation}
\label{3.7}
\Phi_k(z)=\frac{(\kappa_k+1)^{-1}}{\pi i}\int_{L_k}\frac{q_k(t)dt}{t-z}+\frac1{2\pi}\int_{\Gamma_k}\left(g'_k(t)+\frac{2id_k^{-1}q(t)}{\kappa_k+1}\right)\frac{dt}{t-z},
\end{equation} 
$$
\Psi_k(z)=\frac{(\kappa_k+1)^{-1}}{\pi i}\int_{L_k}\left(\frac{\kappa_k\overline{q_k(t)dt}}{t-z}-\frac{\bar{t}q_k(t)dt}{(t-z)^2}\right)
$$
$$
{}+\frac1{2\pi}\int_{\Gamma_k}\left[\left(\overline{g'_k(t)}+\frac{2i\kappa_k \overline{q(t)}}{d_k(\kappa_k+1)}\right)\frac{\overline{dt}}{t-z}-
\left(g'_k(t)+\frac{2id_k^{-1}q(t)}{\kappa_k+1}\right)\frac{\bar{t}dt}{(t-z)^2}\right], 
$$
$$
z\in S_k, \,\,\,\, k=1,\ldots,n.
$$

Observe that the representations (\ref{3.7}) can be also used if the patch $S_k$ covers several holes and is attached both along its own boundary $\Gamma_k$ and all the boundaries of the holes which are covered by $S_k$. In this case the line $L_k$ needs to be understood as the union of all the boundaries of the holes which are covered by $S_k$. 

The same representations (\ref{3.7}) are also valid for the patches $S_k$, $k=n+1,\ldots,n+m$, assuming that $q_k(t)=0$.

Hence, the stressed state of the ``plate-patches" system is described by the complex potentials (\ref{3.4}), (\ref{3.7}) which contain unknown functions $g'(t)$, $t\in L_j$, $j=1,\ldots,n+r$; $q_j(t)$, $t\in L_j$, $j=1,\ldots,n$, and $q(t)$, $g'_k(t)$, $t\in \Gamma_k$, $k=1,\ldots, n+m$. We will look for these functions in the class of functions satisfying the H\"older condition on the corresponding curves $L_j$ and $\Gamma_k$. This choice guarantees the existence of all the principal and the limit values of the integrals of the Cauchy type in the formulas (\ref{3.4}), (\ref{3.7}).

\section{The system of singular integral equations}
To find the unknown functions $g'(t)$, $q_k(t)$, $q(t)$, $g'_k(t)$ we have the conditions (\ref{2.1})-(\ref{2.4}) and (\ref{3.6}). Recall that the condition (\ref{2.5}) has been used to derive representations (\ref{3.7}) and thus is satisfied automatically. Observe also that the condition (\ref{3.6}) appears from the way we extend each patch $S_k$ to the full complex plane (that is, it guarantees that the displacements and the stresses outside of the patch $S_k$ are equal to zero).

From the formulas (\ref{3.1}) and the representations (\ref{3.4}), (\ref{3.7}), satisfying the conditions (\ref{2.1})-(\ref{2.4}) and (\ref{3.6}), we obtain the system of singular integral equations on the closed curves $L_j$, $j=1,\ldots,n+r$, and $\Gamma_k$, $k=1,\ldots,n+m$, with the unknown functions $g'(t)$, $q_k(t)$, $q(t)$, $g'_k(t)$:\\
$$
\frac{\mu_k}{\mu}\left[\kappa\Gamma-\bar{\Gamma}-\bar{\Gamma}'\frac{\overline{dt}}{dt}-\sum_{j=n+1}^{n+r} \left(\frac{\kappa Q_j}{t-z_j}-\frac{\bar{Q}_j}{\bar{t}-\bar{z}_j}+\frac{\overline{dt}}{dt}\left(\frac{\kappa Q_j}{\bar{t}-\bar{z}_j}+\frac{\bar{t}Q_j}{(\bar{t}-\bar{z}_j)^2}\right) \right) \right.
$$
$$
+\frac{i(\kappa+1)}{2}g'(t)+\frac1{2\pi}\int_L\left(\frac{\kappa}{\tau-t}-\frac1{\bar{\tau}-\bar{t}}\frac{\overline{dt}}{dt}\right)g'(\tau)d\tau
$$
$$
+\frac1{2\pi}\int_L\left(-\frac1{\bar{\tau}-\bar{t}}+\frac{\tau-t}{(\bar{\tau}-\bar{t})^2}\frac{\overline{dt}}{dt}\right)\overline{g'(\tau)d\tau}
$$
$$
+\frac{(\kappa+1)^{-1}}{\pi i}\int_{\gamma}\left(\frac{\kappa}{\tau-t}+\frac{\kappa}{\bar{\tau}-\bar{t}}\frac{\overline{dt}}{dt}\right)q(\tau)d\tau
$$
$$
\left.+\frac{(\kappa+1)^{-1}}{\pi i}\int_{\gamma}\left(\frac1{\bar{\tau}-\bar{t}}-\frac{\tau-t}{(\bar{\tau}-\bar{t})^2}\frac{\overline{dt}}{dt}\right)\overline{q(\tau)d\tau}\right]
$$
$$
=\frac{(\kappa_k+1)^{-1}}{\pi i}\int_{L_k}\left(\frac{\kappa_k}{\tau-t}+\frac{\kappa_k}{\bar{\tau}-\bar{t}}\frac{\overline{dt}}{dt}\right)q_k(\tau)d\tau
$$
$$
+\frac{(\kappa_k+1)^{-1}}{\pi i}\int_{L_k}\left(\frac1{\bar{\tau}-\bar{t}}-\frac{\tau-t}{(\bar{\tau}-\bar{t})^2}\frac{\overline{dt}}{dt}\right)\overline{q_k(\tau)d\tau}
$$
$$
+\frac1{2\pi}\int_{\Gamma_k}\left(\frac{\kappa_k}{\tau-t}-\frac1{\bar{\tau}-\bar{t}}\frac{\overline{dt}}{dt}\right)g'_k(\tau)d\tau
$$
$$
+\frac1{2\pi}\int_{\Gamma_k}\left(-\frac1{\bar{\tau}-\bar{t}}+\frac{\tau-t}{(\bar{\tau}-\bar{t})^2}\frac{\overline{dt}}{dt}\right)\overline{g'_k(\tau)d\tau}
$$
$$
+\frac{id_k^{-1}}{\pi(\kappa_k+1)}\int_{\Gamma_k}\left(\frac{\kappa_k}{\tau-t}+\frac{\kappa_k}{\bar{\tau}-\bar{t}}\frac{\overline{dt}}{dt}\right)q(\tau)d\tau
$$
\begin{equation}
\label{3.8}
+\frac{id_k^{-1}}{\pi(\kappa_k+1)}\int_{\Gamma_k}\left(\frac1{\bar{\tau}-\bar{t}}-\frac{\tau-t}{(\bar{\tau}-\bar{t})^2}\frac{\overline{dt}}{dt}\right)\overline{q(\tau)d\tau}, \,\,\,\, t\in L_k,\,\,\, k=1,\ldots,n;
\end{equation}
$$
2d_kq_k(t)+\frac1{2\pi}\int_L\left(\frac1{\tau-t}+\frac1{\bar{\tau}-\bar{t}}\frac{\overline{dt}}{dt}\right)g'(\tau)d\tau
$$
$$
+\frac1{2\pi}\int_L\left(\frac1{\bar{\tau}-\bar{t}}-\frac{\tau-t}{(\bar{\tau}-\bar{t})^2}\frac{\overline{dt}}{dt}\right)\overline{g'(\tau)d\tau}
$$
$$
+\frac{(\kappa+1)^{-1}}{\pi i}\int_{\gamma}\left(\frac1{\tau-t}-\frac{\kappa}{\bar{\tau}-\bar{t}}\frac{\overline{dt}}{dt}\right)q(\tau)d\tau
$$
$$
+\frac{(\kappa+1)^{-1}}{\pi i}\int_{\gamma}\left(-\frac1{\bar{\tau}-\bar{t}}+\frac{\tau-t}{(\bar{\tau}-\bar{t})^2}\frac{\overline{dt}}{dt}\right)\overline{q(\tau)d\tau}
$$
\begin{equation}
=-2 \R \Gamma-\bar{\Gamma}'\frac{\overline{dt}}{dt}+\sum_{j=n+1}^{n+r}\left(2\R\left(\frac{Q_j}{t-z_j}\right)-\frac{\overline{dt}}{dt}\left(\frac{t\bar{Q}_j}{(\bar{t}-\bar{z}_j)^2}+\frac{\kappa Q_j}{\bar{t}-\bar{z}_j}\right) \right),
\label{3.9}
\end{equation}
$$
t\in L_k,\,\,\,\, k=1,\ldots,n;
$$
$$
\frac1{\pi}\frac{|dt|}{dt}\int_{L_k}g'(\tau)d\tau+\frac1{2\pi}\int_L\left(\frac1{\tau-t}+\frac1{\bar{\tau}-\bar{t}}\frac{\overline{dt}}{dt}\right)g'(\tau)d\tau
$$
$$
+\frac1{2\pi}\int_L\left(\frac1{\bar{\tau}-\bar{t}}-\frac{\tau-t}{(\bar{\tau}-\bar{t})^2}\frac{\overline{dt}}{dt}\right)\overline{g'(\tau)d\tau}
$$
$$
+\frac{(\kappa+1)^{-1}}{\pi i}\int_{\gamma}\left(\frac1{\tau-t}-\frac{\kappa}{\bar{\tau}-\bar{t}}\frac{\overline{dt}}{dt}\right)q(\tau)d\tau
$$
$$
+\frac{(\kappa+1)^{-1}}{\pi i}\int_{\gamma}\left(-\frac1{\bar{\tau}-\bar{t}}+\frac{\tau-t}{(\bar{\tau}-\bar{t})^2}\frac{\overline{dt}}{dt}\right)\overline{q(\tau)d\tau}
$$
$$
=-2 \R \Gamma-\bar{\Gamma}'\frac{\overline{dt}}{dt}+\sum_{j=n+1}^{n+r}\left(2\R\left(\frac{Q_j}{t-z_j}\right)-\frac{\overline{dt}}{dt}\left(\frac{t\bar{Q}_j}{(\bar{t}-\bar{z}_j)^2}+\frac{\kappa Q_j}{\bar{t}-\bar{z}_j}\right) \right)+p_k(t),
$$
\begin{equation}
t\in L_k,\,\, k=n+1,\ldots,n+r;
\label{3.10}
\end{equation}
$$
\frac{\mu_k}{\mu}\left[\kappa\Gamma-\bar{\Gamma}-\bar{\Gamma}'\frac{\overline{dt}}{dt}-\sum_{j=n+1}^{n+r} \left(\frac{\kappa Q_j}{t-z_j}-\frac{\bar{Q}_j}{\bar{t}-\bar{z}_j}+\frac{\overline{dt}}{dt}\left(\frac{\kappa Q_j}{\bar{t}-\bar{z}_j}+\frac{\bar{t}Q_j}{(\bar{t}-\bar{z}_j)^2}\right) \right) \right.
$$
$$
+\frac1{\pi}\frac{|dt|}{dt}\int_{\Gamma_k}q(\tau)d\tau+\frac1{2\pi}\int_L\left(\frac{\kappa}{\tau-t}-\frac1{\bar{\tau}-\bar{t}}\frac{\overline{dt}}{dt}\right)g'(\tau)d\tau
$$
$$
+\frac1{2\pi}\int_L\left(-\frac1{\bar{\tau}-\bar{t}}+\frac{\tau-t}{(\bar{\tau}-\bar{t})^2}\frac{\overline{dt}}{dt}\right)\overline{g'(\tau)d\tau}
$$
$$
+\frac{(\kappa+1)^{-1}}{\pi i}\int_{\gamma}\left(\frac{\kappa}{\tau-t}+\frac{\kappa}{\bar{\tau}-\bar{t}}\frac{\overline{dt}}{dt}\right)q(\tau)d\tau
$$
$$
\left.+\frac{(\kappa+1)^{-1}}{\pi i}\int_{\gamma}\left(\frac1{\bar{\tau}-\bar{t}}-\frac{\tau-t}{(\bar{\tau}-\bar{t})^2}\frac{\overline{dt}}{dt}\right)\overline{q(\tau)d\tau}\right]
$$
$$
=\frac{i(\kappa_k+1)}{2}g'_k(t)+\frac{(\kappa_k+1)^{-1}}{\pi i}\int_{L_k}\left(\frac{\kappa_k}{\tau-t}+\frac{\kappa_k}{\bar{\tau}-\bar{t}}\frac{\overline{dt}}{dt}\right)q_k(\tau)d\tau
$$
$$
+\frac{(\kappa_k+1)^{-1}}{\pi i}\int_{L_k}\left(\frac1{\bar{\tau}-\bar{t}}-\frac{\tau-t}{(\bar{\tau}-\bar{t})^2}\frac{\overline{dt}}{dt}\right)\overline{q_k(\tau)d\tau}
$$
$$
+\frac1{2\pi}\int_{\Gamma_k}\left(\frac{\kappa_k}{\tau-t}-\frac1{\bar{\tau}-\bar{t}}\frac{\overline{dt}}{dt}\right)g'_k(\tau)d\tau
$$
$$
+\frac1{2\pi}\int_{\Gamma_k}\left(-\frac1{\bar{\tau}-\bar{t}}+\frac{\tau-t}{(\bar{\tau}-\bar{t})^2}\frac{\overline{dt}}{dt}\right)\overline{g'_k(\tau)d\tau}
$$
$$
+\frac{id_k^{-1}}{\pi(\kappa_k+1)}\int_{\Gamma_k}\left(\frac{\kappa_k}{\tau-t}+\frac{\kappa_k}{\bar{\tau}-\bar{t}}\frac{\overline{dt}}{dt}\right)q(\tau)d\tau
$$
\begin{equation}
+\frac{id_k^{-1}}{\pi(\kappa_k+1)}\int_{\Gamma_k}\left(\frac1{\bar{\tau}-\bar{t}}-\frac{\tau-t}{(\bar{\tau}-\bar{t})^2}\frac{\overline{dt}}{dt}\right)\overline{q(\tau)d\tau}=i(\kappa_k+1)g'_k(t), 
\label{3.11}
\end{equation}
$$
t\in \Gamma_k,\,\,\,\,  k=1,\ldots,n+m;
$$
where $q_k(t)=0$, $k=n+1,\ldots,n+m$, and $d_k=h_k/h$, $k=1,\ldots,n+m$.

Observe, that the additional terms $\frac1{\pi}\frac{|dt|}{dt}\int_{L_k}g'(\tau)d\tau$ and $\frac1{\pi}\frac{|dt|}{dt}\int_{\Gamma_k}q(\tau)d\tau$ have been included into the equations (\ref{3.10}) and (\ref{3.11}) correspondingly. Adding these terms does not change the conditions (\ref{2.1})-(\ref{2.4}) and (\ref{3.6}) and provides that the following conditions are satisfied:
\begin{equation}
\int_{L_k}g'(\tau)d\tau=0,\,\,\, k=n+1,\ldots, n+r,
\label{3.12}
\end{equation} 
\begin{equation}
\int_{\Gamma_k}q(\tau)d\tau=0,\,\,\, k=1,\ldots, n+m.
\label{3.13}
\end{equation} 
From the physical viewpoint  the condition (\ref{3.12}) means that the displacements are single-valued in the plate along the contours $L_k$, $k=n+1,\ldots,n+r$, while the condition  (\ref{3.13}) provides that the resultant force applied to the patch $S_k$ is equal to zero \cite{Zem2007}, \cite{ZemSil2009}. 

It can be shown that the system (\ref{3.8})-(\ref{3.11}) has a unique solution. This is accomplished by reducing the system (\ref{3.8})-(\ref{3.11}) to the system of Fredholm equations and proving the uniqueness of the solution of the mechanical problem. This has been proved in details for the case of the plate with one hole reinforced by one patch in \cite{Zem2007}. The proof for the general case can be obtained using similar methods.

\section{Numerical procedure}

In this section we present an alternative technique for solving the system of singular integral equations (\ref{3.8})-(\ref{3.11}). The method of mechanical quadratures presented to study the problems for cracked plates in \cite{Savruk1981} and adopted for plates with holes reinforced by patches in \cite{Zem2007}, \cite{ZemSil2007}, \cite{ZemSil2009} provides a sufficient accuracy with a reasonable computational time for the systems of similar type in the cases when the number of contours $L_j$ and $\Gamma_k$ is small. However, the results become unsatisfactory both in the sense of accuracy and computational time if a larger number of contours is considered. Another observed disadvantage of the method of mechanical quadratures is that the results obtained using this method tend to deteriorate significantly with the loss of the symmetry in the construction.

To overcome these disadvantages the following approach is proposed and shown to be computationally effective. First, parametrize each curve $L_j$ and $\Gamma_k$ with a parameter $\theta\in[0,2\pi]$. Let $\tau_j=\tau_j(\theta)$ be a parametrization of the curves $L_j$, and $\xi_k=\xi_k(\theta)$ be a parametrization of the curves $\Gamma_k$. Next, approximate the unknown functions $g'(t)$, $q_k(t)$, $t\in L_j$, and $q(t)$, $g'_k(t)$, $t\in\Gamma_k$, by the partial sums of the complex Fourier series by the parameter $\theta$:
$$
g'(t)=\sum_{m=-N}^{N}g^m_{0j}e^{im\theta},\,\,\,\,t\in L_j,\,\,\,\,j=1,\ldots,n+r;
$$
\begin{equation}
q_j(t)=\sum_{m=-N}^{N}q^m_{j}e^{im\theta},\,\,\,\,t\in L_j,\,\,\,\,j=1,\ldots,n;
\label{3.14}
\end{equation}
$$
q(t)=\sum_{m=-N}^{N}q^m_{0j}e^{im\theta},\,\,\,\,t\in \Gamma_j,\,\,\,\,j=1,\ldots,n+m;
$$ 
$$
g'_j(t)=\sum_{m=-N}^{N}g^m_{j}e^{im\theta},\,\,\,\,t\in \Gamma_j,\,\,\,\,j=1,\ldots,n+m,
$$
where $g^m_{0j}$, $q^m_{j}$, $q^m_{0j}$, $g^m_{j}$ are the unknown complex coefficients of the partial sums (\ref{3.14}) which need to be found. 

Substituting the formulas (\ref{3.14}) into the system (\ref{3.8})-(\ref{3.11}) leads to the following equations:
$$
\sum_{k=1}^{n+r}\left(\sum_{m=-N}^N g^m_{0k}A^m_{j,k}(\theta)+\sum_{m=-N}^N \overline{g^m_{0k}}B^m_{j,k}(\theta)\right)
$$
$$
+\sum_{k=1}^{n+m}\left(\sum_{m=-N}^N q^m_{0k}A^m_{j,k+n+r}(\theta)+\sum_{m=-N}^N \overline{q^m_{0k}}B^m_{j,k+n+r}(\theta)\right)
$$
\begin{equation}
+\sum_{k=1}^{n}\left(\sum_{m=-N}^N q^m_{k}A^m_{j,k+2n+m+r}(\theta)+\sum_{m=-N}^N \overline{q^m_{k}}B^m_{j,k+2n+m+r}(\theta)\right)
\label{3.15}
\end{equation}
$$
+\sum_{k=1}^{n+m}\left(\sum_{m=-N}^N g^m_{k}A^m_{j,k+3n+m+r}(\theta)+\sum_{m=-N}^N \overline{g^m_{k}}B^m_{j,k+3n+m+r}(\theta)\right)=f_j(\theta),
$$
$$
\theta\in[0,2\pi],\,\,\,j=1,\ldots,4n+2m+r,
$$
where $N$ is the degree of the trigonometric polynomials (\ref{3.14}), and the functions $A^m_{j,k}(\theta)$, $B^m_{j,k}(\theta)$ and $f_j(\theta)$ can be obtained from the system (\ref{3.8})-(\ref{3.11}). The integrals in these functions can be computed numerically by using any appropriate method.

The equations (\ref{3.15}), in general, can not be satisfied for all the values of the parameter $\theta$ due to the approximation of the unknown functions by the truncated Fourier series (\ref{3.14}). Instead, we assume that the equations (\ref{3.15}) are satisfied only in the finite number of points $\theta_s=\pi (2s-1)/(2N+1)$, $s=1,\ldots,2N+1$. This transforms (\ref{3.15}) into the system of linear algebraic equations with real coefficients of the order $2(2N+1)(4n+2m+r)$. Thus, the approximate solution of the system of singular integral equations (\ref{3.8})-(\ref{3.11}) can be obtained by solving the system of linear algebraic equations.

\begin{figure}
	\centering
		\scalebox{0.5}{\includegraphics{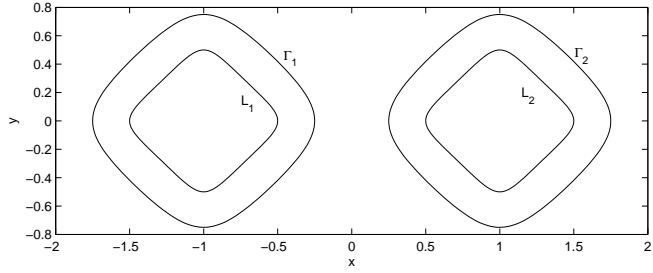}}
	\caption{A plate with two square holes reinforced by two square patches.}
	\label{fig2}
\end{figure}

While the precise investigation of the convergence of the solution of the resulting system of linear algebraic equations to the solution of the system of singular integral equations (\ref{3.8})-(\ref{3.11}) is outside of the scope of this paper, the numerical results obtained by using this method and also comparison with the numerical results obtained by the method of mechanical quadratures in \cite{ZemSil2009} suggests that even relatively small number of terms $N$ ($N\leq 20$ in most cases) is sufficient to obtain a satisfactory convergence, thus, reducing the problem to the solution of a relatively small system of linear algebraic equations. 

\begin{figure}
	\centering
		\scalebox{0.5}{\includegraphics{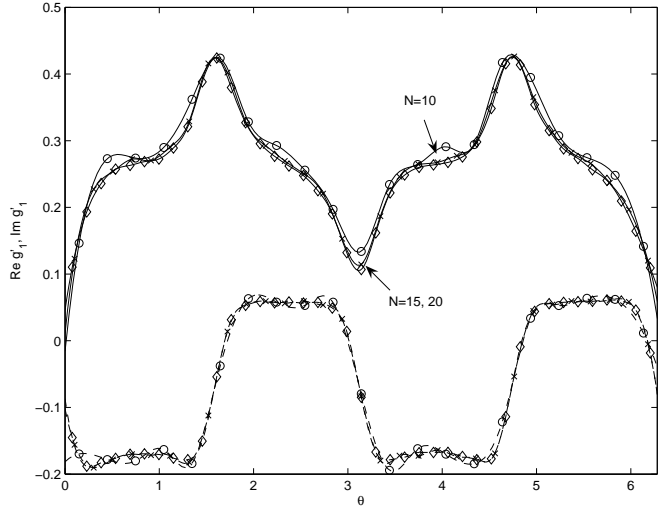}}		
	\caption{The comparison of the approximations of the real and imaginary parts of the function $g'_1(\tau_1(\theta))$ for $N=10$, $N=15$ and $N=20$.}
	\label{fig3}
\end{figure}

\begin{figure}
	\centering
		\scalebox{0.5}{\includegraphics{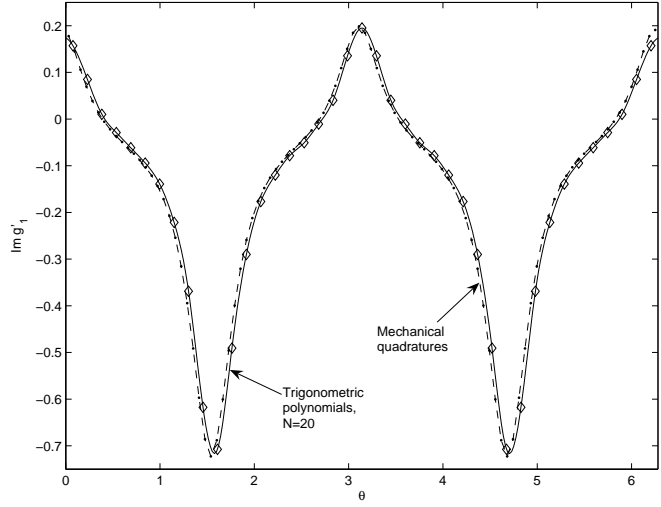}}		
	\caption{The comparison of the approximation of the function $\I g'_1(\tau_1(\theta))$ computed by the method of trigonometric polynomials and the method of mechanical quadratures \cite{ZemSil2009}.}
	\label{fig8}
\end{figure}

For instance, consider a plate with two holes $L_1$ and $L_2$ in the shape of the squares with rounded corners given by the equations $\tau_{1,2}(\theta)=0.45(e^{i\theta}+e^{-3i\theta}/9)\mp 1$ which is reinforced by two patches in the shape of the squares with rounded corners with the boundaries $\Gamma_1$ and $\Gamma_2$ given by the equations $\xi_{1,2}(\theta)=0.7(e^{i\theta}+e^{-3i\theta}/14)\mp 1$ (fig. \ref{fig2}). Assume that the patches attached to the plate along both contours $L_1$, $\Gamma_1$ and $L_2$, $\Gamma_2$ correspondingly. The plate is loaded at infinity by the tensile stress $\sigma_1^{\infty}$ acting in the direction constituting the angle $\alpha=\pi/4$ with the positive direction of $x$-axis, and $\sigma_2^{\infty}=0$.  The mechanical parameters of the plate and of the patches are $\mu=60$, $\nu=0.4$ and $\mu_1=\mu_2=40$, $\nu_1=\nu_2=0.3$ correspondingly. The thicknesses of the patches and of the plate are equal $h=h_1=h_2$. The fig. \ref{fig3} shows the approximations of the real (solid line) and imaginary (dashed line) parts of the unknown function $g'_1(\tau_1(\theta))$ plotted for different numbers $N$ of terms of the trigonometric polynomials (\ref{3.14}): $N=10$, $N=15$ and $N=20$. It can be seen that a good approximation is obtained already for $N=15$. Similar results have been observed in other considered examples. The comparison of the results with the method of mechanical quadratures \cite{ZemSil2009} is presented on the fig. \ref{fig8}. The computation is made for the same mechanical and geometric parameters of the plate and the patches as in fig. \ref{fig2}, and $\sigma_1^{\infty}=1$, $\sigma_2^{\infty}=0$, $\alpha=0$. The patches are joined with the plate only along their boundaries $\Gamma_1$ and $\Gamma_2$. The solid-diamond line corresponds to the results obtained by the method of trigonometric polynomials with $N=20$ presented here and the dash-dot line corresponds to the method of mechanical quadratures with $100$ points on each of the lines \cite{ZemSil2009}. It can be seen from the fig. \ref{fig8} that a good comparison of the results is observed. At the same time it needs to be stressed that the method of mechanical quadratures tends to lose its efficiency for a larger number of contours and that is why alternative numerical techniques should be used. 

\begin{figure}
	\centering
		\scalebox{0.6}{\includegraphics{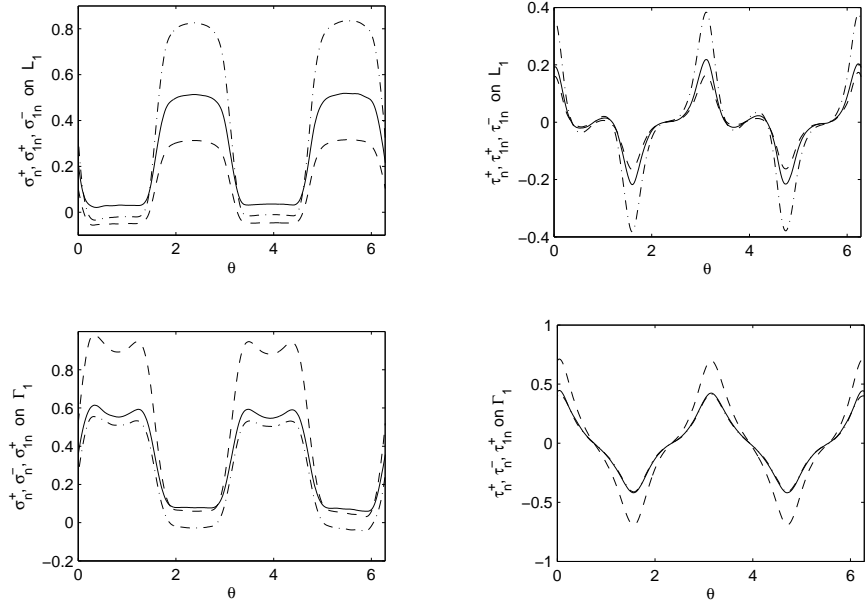}}
	\caption{The stresses on the lines $L_1$ and $\Gamma_1$ in the case of two square holes and two square patches.}
	\label{fig4}
\end{figure} 

\begin{figure}
	\centering
		\scalebox{0.6}{\includegraphics{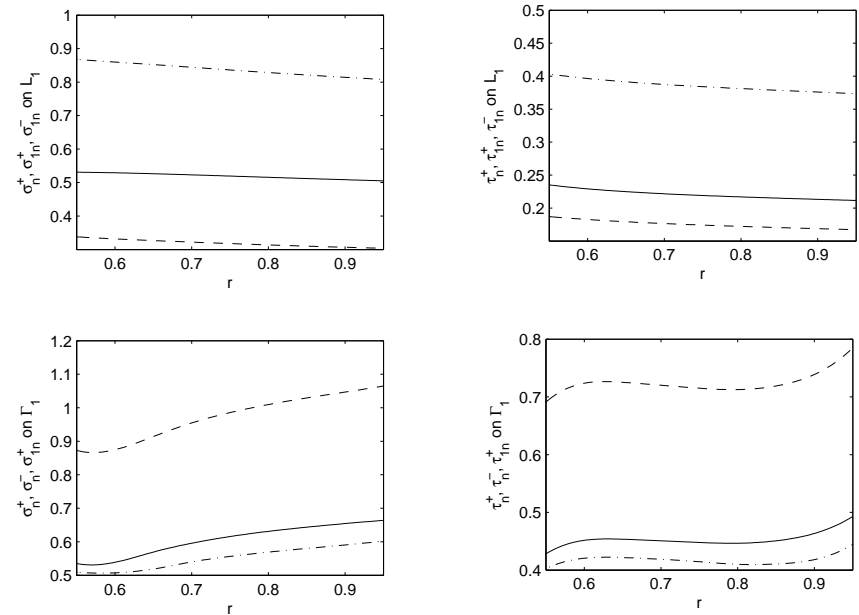}}
	\caption{The dependence of the stresses on the lines $L_1$ and $\Gamma_1$ on the size of the patch.}
	\label{fig7}
\end{figure} 

The fig. \ref{fig4} shows the graphs of the tensile and the shear stresses $\sigma_n$ and $\tau_n$ on the boundaries of the left hole $L_1$ and the left patch $\Gamma_1$ for the construction shown on the fig. \ref{fig2}. The plate is loaded at infinity by the tensile stress $\sigma_1^{\infty}=1$ acting in the direction constituting the angle $\alpha=\pi/4$ with the positive direction of the real axis, and $\sigma_2^{\infty}=0$. The mechanical parameters of the plate and the patches are the same as in the previous example. Here and further, on the hole boundary $L_j$, the graphs given by the solid line represent the stresses in the plate, the graphs given by the dash line represent the stresses in the patch from the outside of the hole boundary $L_j$, and the graphs given by the dash-dot line represent the stresses in the patch from the inside of the hole boundary $L_j$. Similarly, on the patch boundary $\Gamma_j$, the graphs given by the solid line represent the stresses in the plate from the inside of $\Gamma_j$, the graphs given by the dash line show the stresses in the plate from the outside of $\Gamma_j$, and finally, the graphs given by the dash-dot line represent stresses in the patch.

\begin{figure}
	\centering
		\scalebox{0.6}{\includegraphics{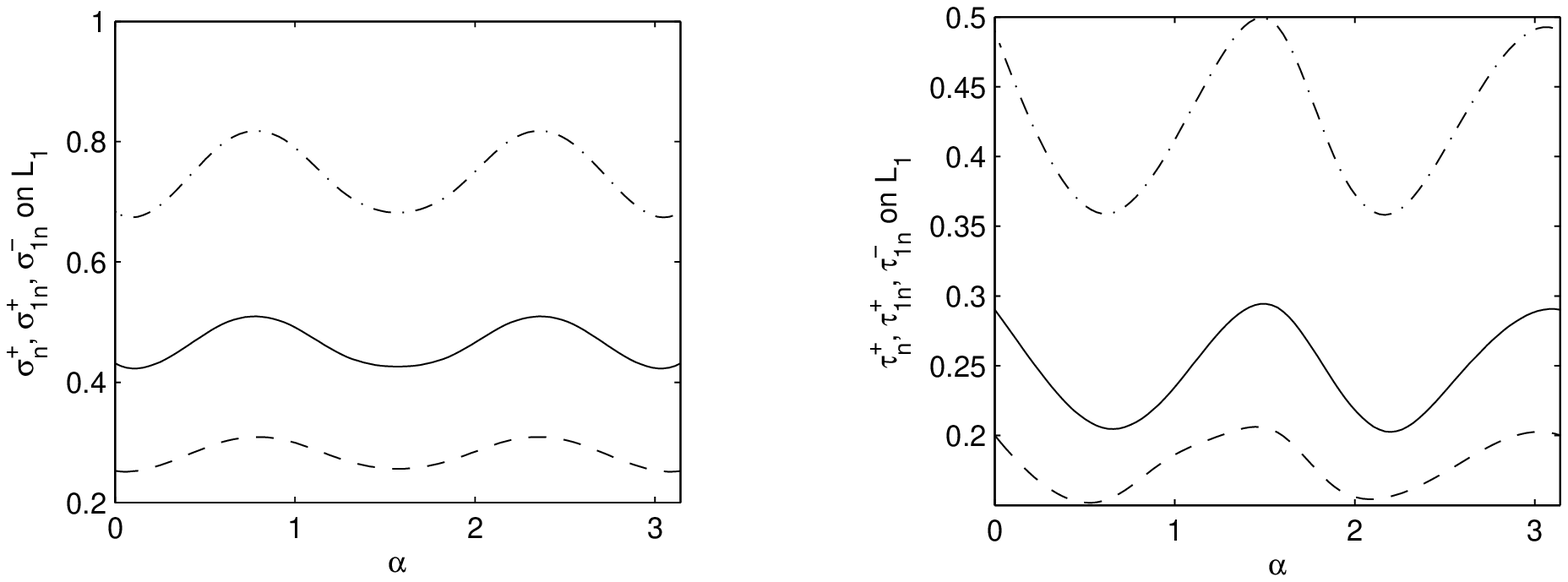}}
	\caption{The dependence of the stresses on the line $L_1$ on the angle $\alpha$.}
	\label{fig5}
\end{figure}   

\begin{figure}
	\centering
		\scalebox{0.6}{\includegraphics{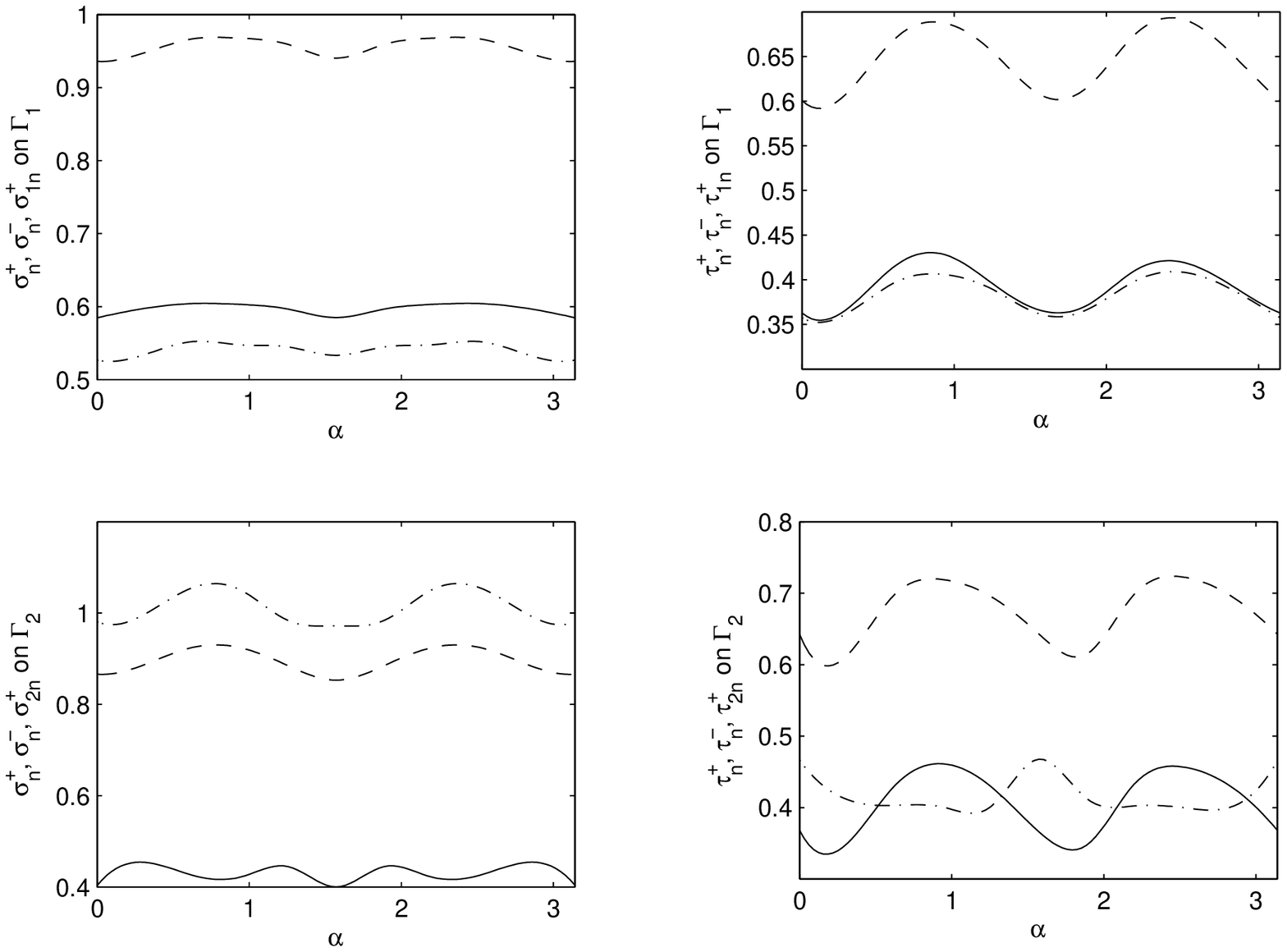}}
	\caption{The dependence of the stresses on the lines $\Gamma_1$ and $\Gamma_2$ on the angle $\alpha$.}
	\label{fig6}
\end{figure}  

The graphs of the dependence of the maximums of the stresses $\sigma_n$ and $\tau_n$ on the lines $L_1$ and $\Gamma_1$ on the patch size characterized by the parameter $r$ are shown on the fig. \ref{fig7}. The parameter $r$ here is the radius of the circle circumscribed around each square patch. The proportions of the patch stay the same as in the previous example and only the size of the patch is changing. The graphs on the lines $L_2$ and $\Gamma_2$ are the same due to the symmetry of the construction. It can be seen from the graphs that while the maximums of the stresses on the hole boundary $L_1$ decrease as the size of the patch increases, the maximums on the patch boundary $\Gamma_1$ increase with $r$. 

\begin{figure}
	\centering
		\scalebox{0.5}{\includegraphics{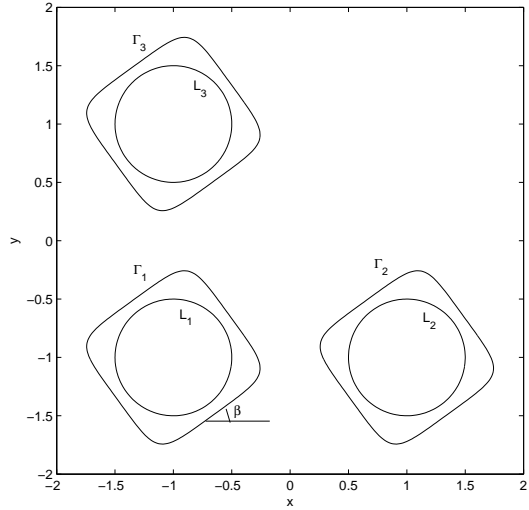}}
	\caption{A plate with three circular holes reinforced by three square patches.}
	\label{fig9}
\end{figure}

Consider the same construction as on the fig. \ref{fig2} where the left patch $S_1$ is attached to the plate $S$ along both lines $L_1$ and $\Gamma_1$ while the right patch $S_2$ is attached only along its boundary $\Gamma_2$. All other mechanical and geometrical parameters are the same as in the previous example. The plate $S$ is loaded at infinity with the tensile stress $\sigma_1^{\infty}=1$ acting in the direction constituting the angle $\alpha$ with the positive direction of the real axis. The dependence of the maximums of the stresses $\sigma_n$ and $\tau_n$ on the line $L_1$ on the direction of the loading given by the angle $\alpha$ is shown on the fig. \ref{fig5}. Observe that the stresses $\sigma_n$ and $\tau_n$ are equal to zero on the line $L_2$. The graphs of the stresses on the lines $\Gamma_1$ and $\Gamma_2$ are shown on the fig. \ref{fig6}. It can be seen that the direction of loading makes a significant impact on the maximums of the stresses and, hence, on the possibility of the failure of the construction. Thus, each potential repair should be made with a consideration for the direction of the loading which is applied to the construction. 

\begin{figure}
	\centering
		\scalebox{0.6}{\includegraphics{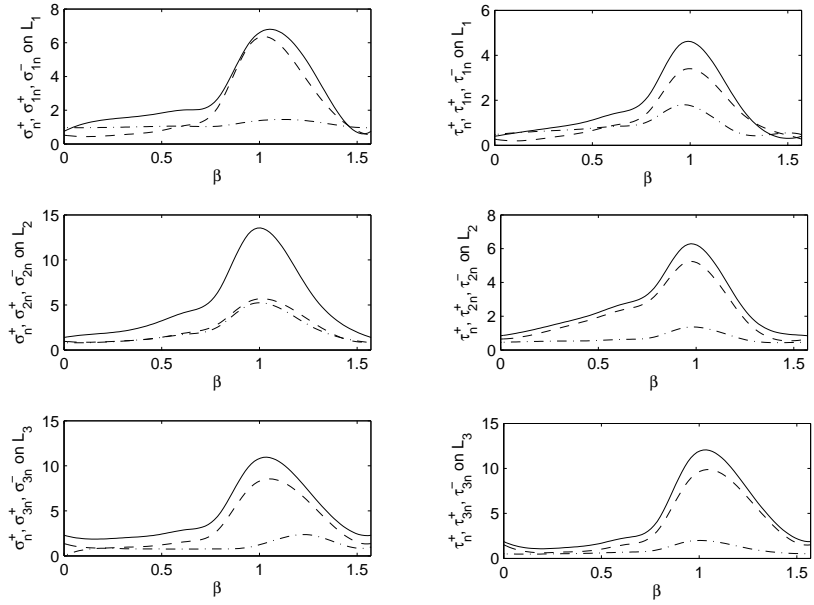}}
	\caption{The dependence of the stresses on the lines $L_1$, $L_2$ and $L_3$ on the angle $\beta$.}
	\label{fig10}
\end{figure}   

\begin{figure}
	\centering
		\scalebox{0.6}{\includegraphics{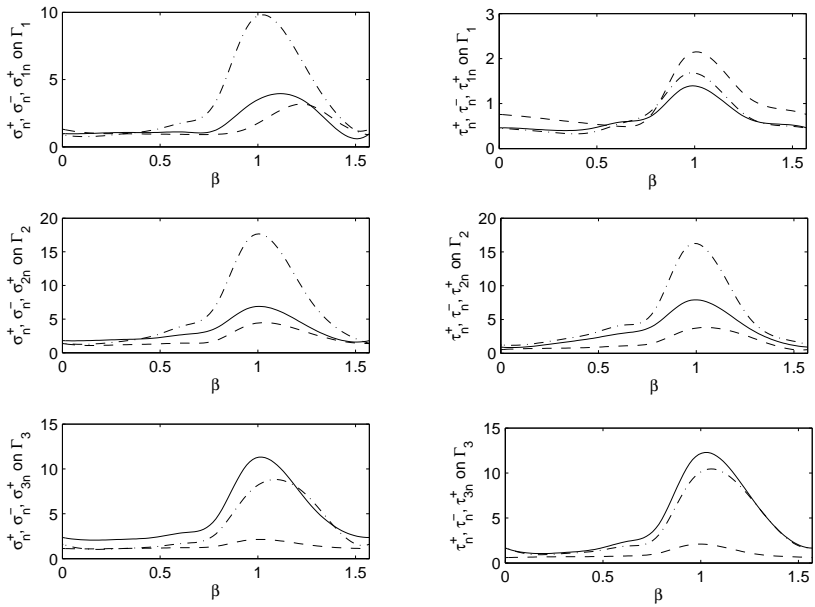}}
	\caption{The dependence of the stresses on the lines $\Gamma_1$, $\Gamma_2$ and $\Gamma_3$ on the angle $\beta$.}
	\label{fig11}
\end{figure}  

The fig. \ref{fig9} shows a plate with three circular holes each of the radius $R=0.5$ and reinforced with three square patches $S_1$, $S_2$, $S_3$ with the rounded corners. The centers of the holes and the corresponding patches are located at the points $z_1=-1-i$, $z_2=1-i$ and $z_3=-1+i$. The boundaries of the patches are given by the equations $\xi_j(\theta)=0.675 \cdot e^{i(\beta-\pi/4)}(e^{i\theta}+e^{-3i\theta}/9)+ z_j$, $j=1,2,3$. The patches are attached to the plate both along their boundaries and the boundaries of the holes they cover.  The angle $\beta$ here represents the orientation of the patches (fig. \ref{fig9}). The plate is loaded at infinity with a normal stress $\sigma_1^{\infty}=1$ acting in the direction of the real axis ($\alpha=0$). The mechanical parameters of the plate and of the patches are $\mu=60$, $\nu=0.4$ and $\mu_1=\mu_2=\mu_3=40$, $\nu_1=\nu_2=\nu_3=0.3$, $h=h_1=h_2=h_3$ correspondingly. 

The graphs of the dependence of the maximal stresses $\sigma_n$ and $\tau_n$ on the lines $L_1$, $L_2$ and $L_3$ on the orientation of the patches $\beta$ are shown on the fig. \ref{fig10}. The corresponding graphs on the lines $\Gamma_1$, $\Gamma_2$ and $\Gamma_3$ are shown on the fig. \ref{fig11}. It can be seen that the maximums of the stresses in the plate and patches depend strongly on the orientation of the patches with respect to the loading applied. The minimums of the stresses are achieved at $\beta\approx 0^{\circ}$, and the maximal stresses achieved at $\beta\approx 60^{\circ}$. There is a noticeable difference in the magnitude of the stresses acting in the construction depending on the orientation of the patches, which means that an improper repair can increase the stresses even compared to the unreinforced case and can potentially lead to a failure of the construction. Thus each potential repair should be carefully evaluated in order to reduce the maximal stresses in the construction. 

The numerical method proposed here for the solution of the systems of singular integral equations of the type (\ref{3.8})-(\ref{3.11}) can be easily extended to more general constructions. It also has been observed to provide more stable numerical results for a wider range of geometries of the holes and patches compared to the method of the mechanical quadratures \cite{Zem2007}, \cite{ZemSil2007}, \cite{ZemSil2009}, to be less dependent on the symmetry of the construction and to be more computationally efficient.

\section{Conclusions}
In this paper we extended a recently suggested technique for solving the problem of an infinite elastic plate containing several holes of arbitrary shapes reinforced by multiple patches of arbitrary shapes to the case when each patch can be attached to the plate by one of the two methods: either only along the boundary of the patch or both along the boundary of the patch and the boundaries of the holes covered by the patch. The problem is reduced to the system of singular integral equations. A new and more efficient numerical approach is proposed for numerical solution of this system. The numerical examples for some particular cases of the plate and patch geometries are given. Obvious future developments of the method are the implementation of finite boundaries of the plate and inhomogeneously imperfect interfaces and also the optimization of the patch repair in accordance with a given plate geometry and an applied loading.

\section*{Acknowledgements}
This work was
supported by Award No. KUS-C1-016-04, made by King Abdullah University of Science and Technology (KAUST).


\vspace{.1in}
 
 
\begin{thebibliography}{99}

\bibitem{BakJon1988} {\sc A.A. Baker, R. Jones,} {\em Bonded repair of aircraft structures}, Martinus Nijhoff Publishers, Dordrecht, 1988. 

\bibitem{EngZakBec2001} {\sc H. Engels, D. Zakharov, W. Becker,} {\em The plane problem of an elliptically reinforced circular hole in an anisotropic plate or laminate}, Archive of Applied Mathematics 71 (2001) 601--612. 

\bibitem{ErdArin1972} {\sc F. Erdogan, K. Arin,} {\em A sandwich plate with a part-through and debonding crack}, Engineering Fracture Mechanics 4(3) (1972) 449--458. 

\bibitem{MWCh1975} {\sc R.A. Mitchell, R.M. Wooley, D.J. Chwirut,} {\em Analysis of composite-reinforced cutouts and cracks}, AIAA Journal 13(6) (1975) 744--749.

\bibitem{Mus1963} {\sc N. I. Muskhelishvili,} {\em Some basic problems of the mathematical theory of elasticity; fundamental equations, plane theory of elasticity, torsion, and bending,} Noordhoff International Publishing, Groningen, 1963.

\bibitem{Rose1988} {\sc L.R.F. Rose,} {\em Theoretical analysis of crack patching}, in: Bonded Repair of Aircraft Structures, Martinus Nijhoff Publishers, 1988, pp. 77--105. 

\bibitem{Savruk1981} {\sc M.P. Savruk,} {\em Two-dimensional problems of elasticity for cracked solids,} Naukova dumka, Kiev (in Russian), 1981.

\bibitem{Zem2007} {\sc A.Y. Zemlyanova,} {\em Singular integral equations for a patch repair problem,} International Journal of Solids and Structures 44 (2007) 6860--6877.

\bibitem{ZemSil2007} {\sc A.Y. Zemlyanova, V.V. Silvestrov,} {\em The problem of the reinforcement of a plate with a cutout by a two-dimensional patch,} Journal of Applied Mathematics and Mechanics 71(1) (2007) 40--51.

\bibitem{ZemSil2009} {\sc A.Y. Zemlyanova, V.V. Silvestrov,} {\em Interaction of elastic plate weakened by multiple holes with a set of patches,} Quarterly Journal of Mechanics and Applied Mathematics 62(2) (2009) 201--220. 


\end{thebibliography}
\end{document}